\documentclass[12pt, leqno]{article}
\usepackage{amssymb,amsmath, amscd}

=\oddsidemargin

\font\teneufm=eufm10
\font\seveneufm=eufm7
\font\fiveeufm=eufm5
\newfam\eufmfam
\textfont\eufmfam=\teneufm
\scriptfont\eufmfam=\seveneufm
\scriptscriptfont\eufmfam=\fiveeufm
\def\frak#1{{\fam\eufmfam\relax#1}}
\let\goth\mathfrak
\def\G{\text{\bf  G}}
\def\F{\text{\bf  F}}
\def\bN{\text{\bf  N}}
\def\bG{\text{\bf  G}}
\def\bH{\text{\bf  H}}
\def\bT{\text{\bf  T}}

\def\gg{\goth g}
\def\gs{\goth s}
\def\gh{\goth h}

\def\gl{\goth l}

\def\gp{\goth p}

\def\gq{\goth q}
\def\gz{\goth z}
\def\go{\goth o}
\def\ga{\goth a}

\def\kalg{k \text{\it -alg}}
\def\pp{\mbox{\bf p}}

\def\bAut{\mbox{\bf Aut}}

\def\bGL{\mbox{\bf GL}}

\def\SL{\mbox{\rm SL}}
\def\bSL{\mbox{\bf SL}}
\def\bGL{\mbox{\bf GL}}

\def\beq{\begin{equation}}
\def\eeq{\end{equation}}
\def\bea{\begin{eqnarray}}
\def\eea{\end{eqnarray}}
\def\beas{\begin{eqnarray*}}
\def\eeas{\end{eqnarray*}}

\def\cplus{\hbox{$\supset${\raise1.05pt\hbox{\kern -0.55em
${\scriptscriptstyle +}$}}\ }}

\DeclareMathOperator{\Span}{Span}

\DeclareMathOperator{\Hom}{Hom}
\DeclareMathOperator{\Aut}{Aut}

\DeclareMathOperator{\im}{im}

\DeclareMathOperator{\Id}{Id}

\DeclareMathOperator{\Ad}{Ad} 
 
\DeclareMathOperator{\ad}{\rm{ad}}

\newtheorem{theorem}[equation]{Theorem}

\newtheorem{lemma}[equation]{Lemma}

\newtheorem{corollary}[equation]{Corollary}

\newtheorem{proposition}[equation]{Proposition}

\newtheorem{remark}[equation]{Remark}

\def\Z{\mathbb Z}

\begin{document}

\title{Automorphisms of toroidal Lie superalgebras}
\author{Dimitar Grantcharov and
Arturo Pianzola\thanks{Supported by the NSERC Discovery Grant
program.}}
\date{}

\maketitle

\begin{abstract} We give a detailed description of the algebraic group
$\bAut(\gg)$ of automorphisms of a simple finite dimensional Lie
superalgebra $\gg$ over an algebraically closed field $k$ of
characteristic $0$. We also give a description of the group of
automorphism of the $k$-Lie superalgebra $\gg \otimes_k R$
whenever $R$ is a Noetherian domain with trivial Picard group.
\end{abstract}

\section{Introduction}

In this article we provide a description of the algebraic group
$\bAut(\gg)$ of automorphisms of a finite dimensional simple Lie
superalgebra over an algebraically closed field $k$ of
characteristic $0$. The results herein somehow refine and complement
those of \cite{S} and \cite{GP}. In the Lie algebra case, the group
$\bAut(\gg)$ is a split extension of a finite constant group (the
symmetries of the Dynkin diagram) by a simple group (the adjoint
group, which is also the connected component of the identity of
$\bAut(\gg)$). By contrast, in the super case the analogous
extension is not split, and the connected component of the identity
of $\bAut(\gg)$ need not even be reductive (let alone simple).

Our approach is to view $\gg$ as a module over a Levi subalgebra
$\gg_0^{ss}$ of the even part of $\gg$. We will introduce three
subgroups of $\bAut(\gg)$; denoted by $\bAut(\gg; \gg_0^{ss})$,
$\bAut(\gg, \Pi_0)$ and $\bH$, which help clarify the nature of
$\bAut(\gg)$ and its outer part. These subgroups are interesting
on their own right, and should prove useful in any future
classification of multiloop algebras of $\gg$ via Galois
cohomology (See \cite{P4} and \cite{P5} for the case of affine Lie
algebras, and \cite{GP} for the case of affine Lie superalgebras.
See also \cite{GiPi1} and \cite{GiPi2} for toroidal Lie algebras).

 As another application, we give an explicit description
of the group of automorphisms of the (in general infinite
dimensional) $k$-Lie superalgebra $\gg \otimes_k R$ for a large
class of commutative ring extensions $R/k$. This result generalizes
the simple Lie algebra situation studied in \cite{P1} (which is
considerably easier by comparison). Of particular interest is the
``toroidal case", namely when $R = k[t_1^{\pm 1}, ..., t_n^{\pm
1}]$.

\section{Notation and conventions}
Throughout $k$ will  be an algebraically closed field of
characteristic zero. The category of associative commutative
unital $k$-algebras will be denoted by $k${\it -alg}. If $V$ is a
vector space over $k$, and $R$ an object of $\kalg$, we set
$V_R=V(R) := V \otimes_k R$. For a nilpotent Lie algebra $\ga$, a
finite dimensional $\ga$-module $M$, and  $\lambda \in \frak a^*$,
we denote by $M^\lambda$ the subspace of $M$ on which $a - \lambda(a)$ acts
nilpotently for every $a \in \ga$. We have then $M = \bigoplus_{\lambda \in \ga^*} M^{\lambda}$.
\bigskip

In what follows, $\gg = \gg_{\bar{0}} \oplus \gg_{\bar{1}}$ will
denote a simple finite dimensional Lie superalgebra  over $k$ (see
\cite{K} and \cite{Sch} for details). A {\it Cartan subsuperalgebra}
$\gh = \gh_{\bar{0}} \oplus \gh_{\bar{1}}$ of $\gg$, is by
definition a selfnormalizing nilpotent subsuperalgebra. Then
$\gh_{\bar{0}}$ is a Cartan (in particular nilpotent) subalgebra of
$\gg_{\bar{0}},$ and $\gh_{\bar{1}}$ is the maximal subspace of
$\gg_{\bar{1}}$ on which $\gh_{\bar{0}}$ acts nilpotently (see
Proposition 1 in \cite{PS} for the proof). We denote by $\Delta =
\Delta_{(\gg, \gh)}$ the {\it roots of $\gg$ with respect to $\gh$}.
Thus $\Delta = \{ \alpha \in \gh_{\bar{0}}^*, \; \alpha \neq 0 \; |
\; \gg^{\alpha} \neq 0\}$. For $\bar{\imath} \in \Z/2 \Z$ we set
$\Delta_{\bar{\imath}} = \{ \alpha \in \gh_{\bar{0}}\; | \;
\gg_{\bar{\imath}}^{\alpha} \neq 0\}$. Then $\Delta =
\Delta_{\bar{0}} \cup \Delta_{\bar{1}}$. The {\it root lattice} $\Z
\Delta$ of $(\gg, \gh)$ will be denoted by $Q_{(\gg, \gh)}$.

A linear algebraic group $\bG$ over $k$ (in the sense of
\cite{Borel}) can be thought as a smooth affine algebraic group (in
the sense of \cite{DG}) via its functor of points $\Hom_k(k[\bG] ,
-)$.  We will find both of these points of view useful, and will
henceforth refer to them simply as ``Algebraic Groups" (and trust
that the reader will be able at all times to understand which of
these two viewpoints is being taken).

Let $\Aut_k(\gg)$ be the (abstract) group of automorphisms of
 $\gg$. We point out that by definition, all automorphisms of a Lie superalgebra
preserve the given $\Z/2\Z$-grading. It is clear that $\Aut_k(\gg)$
gives rise to a linear algebraic group over $k$, which we denote by
$\bAut(\gg)$, whose functor of points is given by $\bAut(\gg)(R) =
\Aut_R(\gg_R)$; the automorphisms of the $R$-Lie superalgebra
$\gg_R=\gg \otimes_k R$.

\bigskip

Recall that there are three types of simple finite dimensional Lie
superalgebras.\footnote{These types are not mutually exclusive, and
some overlap is indeed present in small rank.} We use the notation
of \cite{Pen}.
\bigskip

 {\bf Type I:} $\gs \gl (m|n)$, $\gp \gs \gl (r|r)$,  $\go \gs \gp (2| 2n)$, and
 $\gs {\bf p} (l)$ ($m \neq n$, $r \geq 2$, $l \geq 3$). Every Lie superalgebra $\gg$ of type
I comes equipped with a $\Z$-grading  $\gg= \gg_{-1} \oplus
\gg_{0} \oplus \gg_{1}$ with $\gg_{\bar{0}} = \gg_0$ and
$\gg_{\bar{1}} = \gg_{-1} \oplus \gg_1$. In addition,
$\gg_{\bar{0}}$ is reductive and $\gg_{\pm 1}$ are irreducible
$\gg_{\bar{0}}$-modules.
\medskip

{\bf Type II:} $\go \gs \gp (m|2n)$,  $\gp \gs \gq (l)$,  $F(4)$,
$G(3)$, and $D(\alpha)$ ($m \neq 2$, $l \geq 3$). For these Lie
superalgebras $\gg_{\bar{0}}$ is reductive and $\gg_{\bar{1}}$ is
an irreducible $\gg_{\bar{0}}$-module. We set $\gg_{i}:=\gg_{\bar{\imath}}$, $i=0,1$.

\medskip

{\bf Cartan type:} $W(n)$, $S(m)$, $S'(2l)$, $H(r)$ ($n\geq 2$,
$m\geq 3, l\geq 2, r \geq 5$). Every Lie superalgebra $\gg$ of
Cartan type comes equipped with a choice of subspaces $\gg_i$ for
each $i \in \Bbb Z$ (see the Appendix for details). The Lie algebra
$\gg_{\bar{0}}$ is not reductive but admits a $\Z$-grading
$\gg_{\bar{0}} = \gg_{0} \oplus \gg_{2} \oplus ... \oplus \gg_{2r}$
for which $\gg_{0}$ is reductive. The Lie superalgebra $\gg$ itself,
for all $\gg$ except $\gg = S'(2l)$, admits a $\Z$-grading $\gg_{-1}
\oplus \gg_{0} \oplus ... \oplus \gg_{s}$ where $\gg_{-1}$ and
$\gg_s$ are irreducible $\gg_0$-modules. Note that the notation
$\gg_n$ for even $n$ is not ambiguous (i.e., the same space appears
as the degree $n$ component of the $\Z$-gradings of $\gg_{\bar{0}}$
and of $\gg$ mentioned above).

\bigskip

We thus have a $\Z$-grading $\gg_{0} \oplus \gg_{2} \oplus ...
\oplus \gg_{2r}$ of $\gg_{\bar{0}}$ in all cases  ($r=0$ if $\gg$ is
of type I or II). Furthermore, $\gg$ admits a $\Z$-grading $\gg_{-1}
\oplus \gg_{0} \oplus ... \oplus \gg_{s}$ for $\gg$ of type I or
Cartan type except for $\gg  = S'(2l)$. For convenience, all of the
above will be referred to as {\it standard $\Z$-gradings}.

\bigskip

\section{Structure of $\gg$ with respect to $\gg_0^{ss}$}

We set  $\gg_0^{ss} := [\gg_0, \gg_0]$. This is the semisimple part
of the reductive Lie algebra $\gg_0$, and  a Levi subalgebra of
$\gg_{\bar{0}}$. Henceforth we fix a Cartan subsuperalgebra $\gh$ of
$\gg$ for which $\gh_{\gg_0^{ss}}:= \gh_{\bar{0}}\cap \gg_0^{ss}$ is
a Cartan subalgebra of $\gg_0^{ss}$. If $\gg_{\bar{0}}$ is
reductive, then every Cartan subsuperlagebra of $\gg$ has this
property.  For the remaining cases, namely when $\gg$ is of Cartan
type, the choice of $\gh$ is specified in the Appendix. The root
system of $(\gg_0^{ss}, \gh_{\gg_0^{ss}} )$ will be denoted by
$\Delta_{\gg_0^{ss}}$ and the corresponding root lattice by
$Q_{\gg_0^{ss}}$. We fix once and for all a base $\Pi_0$ of
$\Delta_{\gg_0^{ss}}$. Let $p:\gh_{\bar{0}}^* \to \gh_{\gg_0^{ss}}^*
$ be the canonical map (namely the transpose of the inclusion $
\gh_{\gg_0^{ss}} \subset \gh_{\bar{0}})$.

\medskip

\begin{remark} {\rm Assume $V \subseteq \gg$ is an $\gh_{\bar{0}}$-module (under the adjoint action). Then
$V$ can be viewed as an $\gh_{\gg_0^{ss}}$-module as well. The
generalized weight spaces $V^{\lambda}$ are thus defined for both
actions. Note that in the last case, i.e., for $\lambda \in
\gh_{\gg_0^{ss}},$ we have $ V^{\lambda} = \{ x \in V\; | \; [h,x] =
\lambda(x)v, \mbox{ all }h \in \gh_{\gg_0^{ss}} \}. $ Indeed, since
$\gh_{\gg_0^{ss}}$ is a Cartan subalgebra of the semisimple Lie
algebra $\gg_0^{ss}$ it acts semisimply on the $\gg_0^{ss}$-module
$\gg$, hence also on $V$.}
\end{remark}

\medskip
\begin{lemma}\label{roots}

(i) The set of weights of the $(\gg_0^{ss},
\gh_{\gg_0^{ss}})$-module $\gg$ coincides with $p(\Delta)$.
Moreover, $\gg^{p(\alpha)}= \{ x \in \gg \, | \, [h,x] = \alpha (h)x,
\mbox{ for all } h\in \gh_{\gg_0^{ss}}\}$ whenever $\alpha
\in \Delta$.

(ii) The root lattice $Q_{\gg_0^{ss}}$ is a sublattice
of $p(Q_{\gg})$.
\end{lemma}

\noindent {\bf Proof.} (i) Follows from the previous remark.

(ii) Follows from (i), and the fact that $\gg_0^{ss}$ is a
$\gg_0^{ss}$-submodule of $\gg$. \hfill $\square$

\begin{remark}\label{lattices}{\rm We have three different lattices that define
 gradings of $\gg$:

$\bullet$ The root lattice $Q_{\gg}$ of $\gg$;

$\bullet$ The weight lattice $P_{(\gg_0^{ss}, \gh_{\gg_0^{ss}})}$ of
$\gg_0^{ss}$;

$\bullet$ The sublattice $P_{(\gg_0^{ss}, \gh_{\gg_0^{ss}})}(\gg):= \Z
p(\Delta)$ of $P_{(\gg_0^{ss}, \gh_{\gg_0^{ss}})}$ generated by the weights of
the $(\gg_0^{ss},\gh_{\gg_0^{ss}})$-module $\gg$.

We have $Q_{\gg_0^{ss}} \leq P_{(\gg_0^{ss}, \gh_{\gg_0^{ss}})}(\gg)
\leq P_{(\gg_0^{ss}, \gh_{\gg_0^{ss}})}$.  Note that  $P_{(\gg_0^{ss}, \gh_{\gg_0^{ss}})}(\gg)$
might be a proper sublattice of
$P_{(\gg_0^{ss}, \gh_{\gg_0^{ss}})}$: for example, for $\gg = \gp \gs \gq
(n)$, $P_{(\gg_0^{ss}, \gh_{\gg_0^{ss}})}(\gg)=Q_{\gg_0^{ss}}$ and
$\gg_0^{ss} \simeq \gs \gl (n)$.}
\end{remark}

If $V$ is a finite dimensional $k$-space, and $\sigma \in {\rm
GL}_k(V)$, then $\sigma^* \in {\rm GL}_k(V^*)$ will denote the
transpose inverse of $\sigma$. Thus, under the natural pairing
$\langle - \, , - \, \rangle : V^* \times V \to k$ we have
$\langle \sigma^* \alpha, x \rangle = \langle \alpha, \sigma^{-1}x
\rangle$.

\begin{lemma} \label{star}
Assume $\sigma \in \Aut_k (\gg (R))$ stabilizes $
\gh_{\gg_0^{ss}}$ and $\gg_0^{ss}(R)$. Then
$(\sigma_{| \gh_{\gg_0^{ss}}})^*$ stabilizes
$\Delta_{\gg_0^{ss}}$.
\end{lemma}

\noindent {\bf Proof.} A straightforward calculation shows that
$\sigma ((\gg_0^{ss})^{\beta} \otimes R) = (\gg_0^{ss})^{\sigma^*
(\beta)} \otimes R$ for any root $\beta$ in $\Delta_{\gg_0^{ss}}$.
\hfill $\square$

\medskip

\begin{lemma} \label{compatible}
Let $i \in \Z.$ Assume $\alpha_0 \in p(\Delta)$ is such that
$\gg_i^{\alpha_0} \neq 0$. Then there exists a unique $\alpha \in
\Delta$ with the property that $p(\alpha) = \alpha_0$ and
$\gg^{\alpha} \cap \gg_i \neq 0$. Moreover, $\gg_i^{\alpha_0} =
\gg^{\alpha} \cap \gg_i$.
\end{lemma}

\noindent {\bf Proof.} First observe that the spaces $\gg_i$ are
naturally $\gh_{\gg_0^{ss}}$-modules, so the notation
$\gg_i^{\beta}$ is meaningful for all $\beta \in
\gh_{\gg_0^{ss}}$. The inclusion $\gg_i^{p(\alpha)} \supseteq
\gg^{\alpha} \cap \gg_i$ is obvious for every $\alpha \in \Delta$
and $i \in \Z$. To establish the reverse inclusion, we consider
first the case when $\gh_{\bar{0}} \subset \gg_0$. If
$\gh_{\bar{0}} = \gh_{\gg_0^{ss}}$ (i.e. $\gg_0$ is semisimple),
it is clear that $\gg_i^{p(\alpha)} = \gg^{\alpha} \cap \gg_i$ for
every $\alpha \in \Delta$ and $i \in \Z$. Let us assume now that
$\gg_0$ has a nontrivial center $\gz$. In this case we fix $z \in
\gz$ such that $\gz = k z$ and $[z,y]=iy$ whenever $y \in \gg_i$
(see Proposition 1.2.12 in \cite{K}). Considering $\gg_i$ as an
$\gh_{\bar{0}}$-module we have that every nonzero $x_0$ in
$\gg_i^{\alpha_0}$ decomposes as a sum $x_0 = \sum_{\alpha \in
p^{-1}(\alpha_0)\cap \Delta} x_{\alpha}$ for some $x_{\alpha} \in
\gg^{\alpha} \cap \gg_i$ (recall that the action of $
\gh_{\gg_0^{ss}}$ on $\gg$ is semisimple). But since
$[z,x_{\alpha}] = ix_{\alpha}$ we have $\alpha(z) =i$. Therefore
$x_0 \in \gg^{\alpha} \cap \gg_i$, where $\alpha \in \Delta$ is
determined uniquely by
 $\alpha_{|\gh_{\gg_0^{ss}}} = \alpha_0$ and $\alpha(z) = i$. In particular, $\gg_i^{\alpha_0} \subseteq \gg^{\alpha} \cap \gg_i$ and therefore
$\gg_i^{\alpha_0} = \gg^{\alpha} \cap \gg_i$. If $\alpha' \in \Delta$
is another root with the properties described in the lemma, then  for $ 0 \neq x' \in \gg^{\alpha'} \cap \gg_i$ we have $[z,x']=ix'$
and thus $\alpha'(z) = i$. Hence $\alpha'=\alpha$.

It remains to consider the case when $\gh_{\bar{0}} \neq
\gh_{\bar{0}} \cap \gg_0$, which is present for $\gg = H(n)$ only.
Using the explicit description of $\gh_{\bar{0}}$ provided in the
Appendix we see that if $h \in \gh_{\bar{0}}$, then $h = h_0 +
h_n$ for some $h_0 \in \gh_{\gg_0^{ss}}$ and $h_n \in \gg^2 =
\gg_2 \oplus ... \gg_{2r}$ (here $r = \left[
\frac{n-3}{2}\right]$). Let $x$ be a nonzero element in
$\gg_i^{\alpha_0}$ and let $\alpha$ be any root in $\Delta$ such
that $p(\alpha)=\alpha_0$. Then the identities $\alpha(h_n) = 0$
and  $(\ad(h_0) - \alpha(h_0))(x) = 0$ imply that $(\ad(h) -
\alpha(h))^N (x) =(\ad(h_n))^N(x)$ and, in particular, $(\ad(h) -
\alpha(h))^N (x) \in \gg_{2N+i} \oplus \gg_{2N+i+2} \oplus...$
from which it follows that $x \in \gg^{\alpha} \cap \gg_i$ and
thus $\gg_i^{\alpha_0} = \gg^{\alpha} \cap \gg_i$. The uniqueness
of $\alpha$ follows from the fact that $p$ is injective on
$\Delta$ since $\gh_{\bar{0}} = \gh_{\gg_0^{ss}}\oplus \gh^2$ and
$\alpha_{|\gh^2} = 0$ for every $\alpha \in \Delta$.
 \hfill $\square$
\medskip

If $B$ is a basis of $Q_{\gg}$, then $p(B)
\cap Q_{\gg_0^{ss}}$ is not necessarily a basis of
$Q_{\gg_0^{ss}}$. However, this statement is true to some
extent as the following result shows.

\begin{lemma} \label{t2}
Suppose that either $\gg = \gp \gs \gl (2|2)$, or  $\gg \neq
H(2k)$ with $\gg_0^{ss}$ simple.

(i) There exists a basis $B$ of $Q_{\gg}$ for which $p(B) \cap
\Delta_{\gg_0^{ss}}$ is a base of
$\Delta_{\gg_0^{ss}}$.

(ii) Every group homomorphism $\tau$ in $\Hom(Q_{\gg_0^{ss}}, R^{\times})$
can be extended to a group homomorphism $ \widehat{\tau}$ in $
\Hom(Q_{\gg}, R^{\times})$, i.e., $\widehat{\tau} = \tau \circ p$.

(iii) If $\alpha \in B$ is such that $p(\alpha)\in \Delta_{\gg_0^{ss}}$, then
$\gg^{\alpha} \cap \gg_0 = \gg_0^{p(\alpha)} = (\gg_0^{ss})^{p(\alpha)}$.
\end{lemma}
\noindent {\bf Proof.} (i) This statement follows by a case-by-case
verification. We will use the notations and explicit description
of the root systems provided in Appendix A of \cite{Pen}.

\medskip \noindent {\it Case 1: $\gg = \gp \gs \gl (2|2)$.} We set
$B=\{\varepsilon_1 - \varepsilon_2, \varepsilon_{2} - \delta_1,
\delta_1 - \delta_2 \}$.

\medskip  \noindent {\it Case 2: $\gg = \gs \gl (n|1)$.} We set
$B=\{\varepsilon_1 - \varepsilon_2,..., \varepsilon_{n-1} -
\varepsilon_n, \varepsilon_n - \delta \}$.

\medskip \noindent {\it Case 3: $\gg = \go \gs \gp (m|2n)$, $m=1,2$.} In
the case $m=1$ we have that $\Delta =
\Delta_{\gg_0^{ss}}$. For $m=2$ we set $B=\{\varepsilon_1
- \delta_1, \delta_1 - \delta_2,..., \delta_{n-1} - \delta_n, 2
\delta_n \}$.

\medskip \noindent {\it Case 4: $\gg = \gs {\bf p} (n)$.} Set
$B=\{-2\varepsilon_1, \varepsilon_1 - \varepsilon_2,...,
\varepsilon_{n-1} - \varepsilon_n \}$.

\medskip \noindent {\it Case 5: $\gg = \gs \gp \gq  (n)$.} In this case $\Delta =
\Delta_{\gg_0^{ss}}$.

\medskip \noindent {\it Case 6: $\gg = W (n), S(n)$, or $S'(n)$ ($n=2l$ in the last
case)}. Set $B=\{ \varepsilon_1 - \varepsilon_2,...,
\varepsilon_{n-1} - \varepsilon_n, \varepsilon_n \}$.

\medskip \noindent {\it Case 7: $\gg = H(2l+1)$}. We set
$B=\{ \varepsilon_1 - \varepsilon_2,..., \varepsilon_{l-1} -
\varepsilon_l, \varepsilon_l \}$.

\bigskip The second assertion follows directly from (i). It remains to show (iii).
This assertion follows from Lemma \ref{compatible} by verifying
case by case that $\gg^{\alpha} \cap \gg_0 \neq 0$ for $\alpha \in
B$. This verification is trivial if $\gh_{\bar{0}} =
\gh_{\gg_0^{ss}}$ (see the proof of Lemma \ref{compatible}). In
the remaining three cases we proceed as follows.

If $\gg = \gs \gl(n|1)$ then $\gg^{\varepsilon_i-\varepsilon_j}\cap \gg_0 =
\gg_0^{p(\varepsilon_i-\varepsilon_j)}$ is generated by the matrix
$\left(\begin{matrix} E_{ij} & 0 \\ 0 & 0\end{matrix}\right)$, where $E_{ij}$
is the $(i,j)$-th elementary $n\times n$ matrix.

\medskip If $\gg =W(n)$ then $\gg^{\varepsilon_i-\varepsilon_j}\cap \gg_0 =
\gg_0^{p(\varepsilon_i-\varepsilon_j)}$ is generated by the derivations
$D_{\xi_i \xi_j}$.

\medskip If $\gg =H(2l+1)$, then we use the description of  the roots of $\gg$  provided in the Appendix.
In particular, we verify that $\gg^{\varepsilon_i-\varepsilon_j} \cap \gg_0 =
\gg_0^{p(\varepsilon_i-\varepsilon_j)}$
is generated by $D_{\eta_i \eta_j}$, and that $\gg^{\varepsilon_i} \cap \gg_0=
\gg_0^{p(\varepsilon_i)}$ is generated by $D_{\eta_i \eta_{2l+1}}$. \hfill $\square$

\begin{lemma} \label{t3} Suppose that $\gg_0^{ss}$ is not simple and that $\gg \neq
\gp \gs \gl (2|2)$. Then  $\dim \gg^{\alpha} = \dim \gg^{-\alpha}
= 1$ for every $\alpha \in \Delta$. In addition, there exists a set
$\Pi$ of simple roots of $\Delta$ such that the root spaces
$\gg^{\pm \alpha}$, $\alpha \in \Pi$, generate $\gg$.
\end{lemma}

\noindent {\bf Proof.} This follows from the classification
of the contragredient  finite-dimensional Lie superalgebras (see
\S 2.5 and Theorem 3 in \cite{K}). \hfill $\square$

\section{Some  subgroups  of $\bAut(\gg)$}
\setcounter{equation}{0}

In this section we introduce and study an important list of
subgroups of $\bAut(\gg)$ for each simple Lie superalgebra $\gg$. These groups will be used in the next
section for describing $\Aut_k (\gg(R))$.

\subsection{The even superadjoint group $\bG_{\bar{0}}^{\text{\rm sad}}$}

\begin{lemma} \label{ad}
Let $\G$ be the  Chevalley $k$-group of simply connected type
corresponding to $\gg_0^{ss}$. The restriction $\ad_{|\gg_0^{ss}} :
\gg_0^{ss} \to \gg\gl (\gg)$ of the adjoint representation of $\gg$,
lifts uniquely to a morphism $\Ad : \G \to \bAut (\gg)$ of linear
algebraic groups.
\end{lemma}
\noindent {\bf Proof.}  Because $\G$ is simply connected, there
exists a homomorphism $\Ad : \G \to \bGL(\gg)$ of linear algebraic
$k$-groups whose differential is $\ad$. The group $\G$ is generated
by the root subgroups, namely by elements $\exp(z)$, where the
elements $z$ belong to the root spaces of $\gg_0^{ss}$ with respect
to $\gh_{\gg_0^{ss}}$. We have $\Ad (\exp(z)) = \exp(\ad(z))$. Since
the $\exp(\ad(z))$ are automorphisms of the Lie superalgebra $\gg$,
the result follows. \hfill $\square$

\bigskip

The image of $\Ad$ (in the schematic sense) is denoted by
$\bG_{\bar{0}}^{\text{\rm sad}}$: For $R \in k\mbox{-alg}$ and
$\sigma \in \bAut (\gg)(R)$ we have that $\sigma \in
\bG_{\bar{0}}^{\text{\rm sad}}(R)$ if and only if there exists an
$fppf$ extension $\widetilde{R}/R$ and an element $x \in \bG
(\widetilde{R})$ such that $\Ad_{\widetilde{R}}(x) =
\widetilde{\sigma}$, where $\widetilde{\sigma}$  is the image of
$\sigma$ under the map $\bAut (\gg) (R) \to \bAut (\gg)
(\widetilde{R})$. The homomorphism $\Ad$ induces an isomorphism
between the quotient group $\bG/\mbox{ker}(\Ad)$ and
$\bG_{\bar{0}}^{\text{\rm sad}}$. The structure of
$\bG_{\bar{0}}^{\text{\rm sad}}$ is given in Table 1 at the end of
the paper.

\subsection{The diagonal subgroup $\bH$}

 It is  easier to describe $\bH$ via
its functor of points. We have $\bH(R):= \Hom(Q_{\gg},
R^{\times})$. By definition, an element $\lambda$ in $\bH(R)$
fixes $\gh \otimes_k R$ pointwise, and acts on $\gg^{\alpha}
\otimes R$ as a multiplication by $\lambda(\alpha)$. As an
algebraic group, $\bH$ is a split torus (hence connected). In the
case of semisimple Lie algebra, $\bH$ is a Cartan subgroup of adjoint type.

\subsection{The unipotent group $\bN$}

The groups $\bN$ appear only when $\gg$ is of Cartan type. Details
can be found in \S 6.9 of \cite{GP} or in Table 1 at the end. If
$\gg$ is not of Cartan type we set $\bN=1$. For a Cartan type Lie
superalgebra $\gg$ we have
$$
\bN(R):=\{\varphi \in \Aut_R(\gg (R))\; | \; \varphi(x) - x \in
\gg_{2i+2}(R) .... \oplus \gg_{2r}(R), \mbox{ if } x \in
\gg_{2i}\otimes R, i=0,...,r\}.
$$
Alternatively, $\bN$ is the unipotent group that corresponds to
the nilpotent Lie algebra $\gg^2:=\gg_2 \oplus ... \oplus
\gg_{2r}$, whenever $\gg \neq H(2l)$. If $\gg = H(2l)$, then $\bN$
corresponds to $\widetilde{H}(2l)^2:=H(2l)^2\oplus
kD_{\xi_1...\xi_{2l}}$ where $D_f:=\sum_{i=1}^{2l}\frac{\partial f
}{\partial \xi_i}\frac{\partial}{\partial \xi_i}$ whenever  $f$ is
in the Grasmann (super)algebra $\Lambda(\xi_1,...,\xi_{2l})$. Note
that $\gg^2$ is the radical of $\gg_{\bar{0}}$ if $\gg \neq W(n)$,
and is the radical of $[\gg_{\bar{0}}, \gg_{\bar{0}}]$ if $\gg =
W(n)$.

\begin{lemma} \label{nilp}
Let $\gg$ be a Cartan type Lie superalgebra and let $\sigma \in
\bAut(\gg)(R)$. Then $\sigma = \sigma_0 \sigma_n$ where $\sigma_n
\in \bN(R)$ and $\sigma_0$ preserves the standard $\Z$-gradings of
$\gg_{\bar{0}}(R)$ and $\gg(R)$ (the latter if $\gg \neq S'(2k)$).
In particular, $\sigma_0(\gg_0(R)) = \gg_0(R)$.
\end{lemma}
\noindent {\bf Proof.} This follows from the explicit description of
$\bAut(\gg)$ given in \S 6.9 and Table 1 of \cite{GP}. \hfill
$\square$

\subsection{The group $\bSL_2^{\text{\rm out}}$}

 Let  $\gg = \gp \gs \gl (2|2)$.  Recall (see \cite{GP}, \S 6.4, for details)
that we have a closed embedding $\rho : \bSL_2
 \rightarrow \bAut(\gg)$. The image of $\rho$ will be denoted by $\bSL_2^{\text{\rm out}}$.
 We recall for future use the explicit
 nature of $\rho$.

Let $V_2$ be the standard $\gg \gl_2$-module. We have that $\gg =
\gg_{-1} \oplus \gg_{0} \oplus \gg_1$, where $\gg_{0} =
\gg_{\bar{0}} \simeq \gs \gl (2) \oplus \gs \gl (2)$, $\gg_{-1}
\simeq  V_2^* \otimes V_2$, and $\gg_1 \simeq V_2 \otimes V_2^*$.
The $\gg_0$-modules $\gg_1 $ and $\gg_{-1}$ are isomorphic and an
explicit isomorphism $\phi : \gg_1 \to \gg_{-1}$ is determined by
the linear transformation $\psi : {\cal M}_2 \to {\cal M}_2$, where
$\psi(E) := -JE^tJ^{-1}$ for $J = \left(
\begin{array}{cc} 0&1\\-1&0
\end{array} \right)$ (${\cal M}_2$ is the set of $2\times 2$
matrices with entries in $k$). Then $\rho$ is given by
$$
\left( \begin{array}{cc} A & B\\ C & D
\end{array} \right) \mapsto \left( \begin{array}{cc} A &\alpha B +
\beta \psi(C)\\ \gamma \psi(B) + \delta C& D \end{array} \right),
$$
where  $\left(
\begin{array}{cc} \alpha&\beta\\ \gamma&\delta
\end{array} \right) \in \SL_2(k)$.

\subsection{The group $\bAut(\gg ; \gg_0^{ss})$} \label{ptwise}

 For a given Lie subalgebra $\gs$ of $\gg_0^{ss}$, we let
$\bAut(\gg; \gs)$ be the subgroup of $\bAut(\gg)$ consisting of
those automorphisms that fix all elements of $\gs$. Clearly
$\bAut(\gg; \gs)$ is a closed subgroup of $\bAut(\gg)$, hence a
linear algebraic group. Its functor of points is given by
$\bAut(\gg; \gs)(R) = \{\varphi \in \Aut_R (\gg(R))\; | \;
\varphi_{| \gs(R)} = \Id \}$. Our main interest is the case when
$\gs = \gg_0^{ss}$.

 If $\gg = \bigoplus_{n \in \Z} \gg_n$ is
$\Z$-graded, then each $\lambda \in k^{\times}$ defines an
automorphism $\delta_{\lambda} \in \Aut_k (\gg)$ via
$\delta_{{\lambda}|{\gg_n}} := \lambda^n \Id$. This yields a closed
embedding  $\delta : \bG_m \to \bAut(\gg)$. Along similar lines, if
$\gg$ is $\Z/m\Z$-graded, we have a closed embedding $\delta :
\boldsymbol{\mu}_m \to \bAut(\gg)$ (this presupposes a choice of
primitive $m$th root of unity in $k$.) It is clear that for the
standard $\Z$ or $\Z/2\Z$ gradings of $\gg$ the resulting closed
subgroups $\delta(\bG_m)$ or $\delta(\boldsymbol{\mu}_2)$ of
$\bAut(\gg)$ lie inside $\bAut(\gg; \gg_0^{ss})$.

 If $\gg = H(2l)$, we have an extra  ``additive" group of
automorphisms of $\gg$ which we now describe. For $a \in k$ we
define an automorphism $B_{a}$ of $\Lambda(\xi_1,...,\xi_{2l})$,
by $B_{a}(\xi_i):=\xi_i + a \frac{\partial}{\partial
\xi_i}(\xi_1...\xi_{2l})$. This automorphism lifts to an
automorphism $\beta_{a}$ of $\gg$ in the standard way:
$\beta_a(D):= B_a  D B_a^{-1}$ for $D \in \gg$. Since
$B_{a_1}B_{a_2} = B_{a_1 + a_2}$, we have a closed embedding
$\beta: \bG_a \to \bAut(\gg)$. Furthermore, if $\Delta_{\lambda}$
is the automorphism of $\Lambda(\xi_1,...,\xi_{2l})$ corresponding
to $\delta_{\lambda}$, then one easily checks that
$\Delta_{\lambda}(\xi_i) = \lambda \xi_i$, and therefore
$$
B_a (\Delta_{\lambda} (\xi_i)) = \Delta_{\lambda}(B_a (\xi_i)) =
\lambda \xi_i + a \lambda^{2l-1} \frac{\partial}{\partial
\xi_i}(\xi_1...\xi_{2l}).
$$
The latter identity easily implies that $\delta$ and $\beta$
commute, and that the resulting homomorphism $\delta \times \beta:
\bG_m \times \bG_a \to \bAut(\gg; \gg_0^{ss})$ is injective.

\begin{lemma} \label{cart}
Let $\gg$ be a Cartan type Lie superalgebra. The algebraic group
$\mbox{\bf N} \cap  \bAut(\gg; \gg_0^{ss})$ is trivial if
$\gg \neq H(2l)$. If $\gg = H(2l)$, then $\beta: \bG_a \to
\mbox{\bf N} \cap  \bAut(\gg; \gg_0^{ss})$ is an
isomorphism.
\end{lemma}
\noindent {\bf Proof.}  Let $\varphi \in \mbox{\bf N} \cap
\bAut(\gg; \gg_0^{ss}),$ and let $\Phi$ be an automorphism of
$\Lambda(n)$ for which $\varphi(D)= \Phi  D   \Phi^{-1}$. Then
$D(\Phi(\xi_i))= \Phi(D(\xi_i))$ for every $D \in \gg_0^{ss}$ and
$i =1,...,n$. Let $\Phi(\xi_i) = \xi_i + f_i$ with $\deg f_i \geq
3$ or $f_i=0$.

Consider first the case $\gg \neq H(n)$. Then $\xi_j
\frac{\partial}{\partial \xi_i} \in \gg_0^{ss}$ whenever
$1\leq i \neq j \leq n$ which leads to the following chain of
identities
$$
\xi_j + \xi_j\frac{\partial f_i}{\partial \xi_i} = \xi_j
\frac{\partial}{\partial \xi_i} (\xi_i + f_i) = \Phi(\xi_j
\frac{\partial}{\partial \xi_i}(\xi_i)) = \xi_j + f_j.
$$
Therefore $f_j = \xi_j \frac{\partial f_i}{\partial \xi_i}$, and
in particular $\frac{\partial f_j}{\partial \xi_i}=0$ for every $i
\neq j$ The latter implies that $f_j$ is a scalar multiple of
$\xi_j$. Thus $f_j=0$, and therefore $\Phi = \Id$.

Let now $\gg = H(n)$. Then $\gg_0^{ss}$ is generated by
the elements $D_{\xi_i \xi_j} = \xi_j \frac{\partial}{\partial
\xi_i} - \xi_i \frac{\partial}{\partial \xi_j}$, $i\neq j$. As
before, we find
$$
\xi_j + \xi_j \frac{\partial f_i}{\partial \xi_i} - \xi_i
\frac{\partial f_i}{\partial \xi_j} = \left( \xi_j \frac{\partial
}{\partial \xi_i} - \xi_i \frac{\partial }{\partial
\xi_j}\right)(\xi_i + f_i) = \Phi \left( \xi_j \frac{\partial
}{\partial \xi_i} - \xi_i \frac{\partial }{\partial
\xi_j}\right)(\xi_i) = \xi_j +f_j.
$$
Therefore $f_j= D_{\xi_i\xi_j}f_i$. We now make use of the
explicit description of $\bN$ given in \cite{GP}\footnote{More
precisely the proof of  Lemma 6.2 (ii). The key observation is
that the map $f \mapsto D_f$ defines an isomorphism from
$\Lambda^{2i+2} (V)= \Lambda^{2i+2} (\xi_1,...,\xi_n)  $ to the
subspace of homomorphisms in $\Hom (V, \Lambda^{2i+1} (V))$ which
leave invariant the form $\omega_n$ defining  $H(n)$. We then sum
over all $i \geq 0$, and see that every element of $\bN_{H(n)}(k)
\simeq \oplus_{i \geq 0} \Lambda^{2i} (V)$ corresponds to an
automorphism of $\Aut \Lambda(\xi_1,...,\xi_n)$ of the form $f
\mapsto f + \sum_{i\geq 2} D_{F_{2i}}f$.}
 to conclude that
$$
f_i = D_{F_4}(\xi_i)+ D_{F_6}(\xi_i) +... = \frac{\partial F_4
}{\partial \xi_i} + \frac{\partial F_6}{\partial \xi_i}+...
$$
for some $F_i \in \Lambda (\xi_1,...,\xi_n)$ with $\deg F_i = i$.
In particular $\frac{\partial f_i}{\partial \xi_i} = 0$, and thus
$f_j = - \xi_i\frac{\partial f_i}{\partial \xi_j}$ for every $i
\neq j$. Since $f_j$ is a multiple of $\xi_i$ for every $i \neq
j$, we find that $f_j = c_j \xi_1...\xi_{j-1}\xi_{j+1}...\xi_n$
for some constants $c_j$ in $k$. Because of the parity preserving
nature of $\Phi$, we conclude that $c_j=0$ if $n$ is odd.

 Suppose now $\gg = H(2l)$. We first rewrite $f_j$ as $f_j = c_j
\xi_1...\xi_{j-1}\xi_{j+1}...\xi_n = (-1)^{j-1}c_j \frac{\partial
F}{\partial \xi_j}$ where $F:=\xi_1...\xi_{2l}$. Put $a_j :=
(-1)^{j-1} c_j$. Using that $f_j = a_j \frac{\partial F}{\partial
\xi_j}$, $f_i = a_i \frac{\partial F}{\partial \xi_i}$, and $f_j =
- \xi_i\frac{\partial f_i}{\partial \xi_j}$, we verify that
$$
a_j \frac{\partial F}{\partial \xi_j} = f_j =  -
\xi_i\frac{\partial f_i}{\partial \xi_j} = - \xi_i
\frac{\partial}{\partial \xi_j} \left( a_i \frac{\partial
F}{\partial \xi_i}\right) = - a_i \xi_i \frac{\partial}{\partial
\xi_j} \frac{\partial F}{\partial \xi_i} =  a_i
\frac{\partial}{\partial \xi_j} \left( \xi_i\frac{\partial
F}{\partial \xi_i} \right) = a_i \frac{\partial F}{\partial \xi_j}
$$
(note that $\xi_j \frac{\partial F}{\partial \xi_j} = F$ for every
$j$, and that $\xi_i \frac{\partial }{\partial \xi_j} = -
\frac{\partial }{\partial \xi_j} \xi_i$ for every $i$ and $j$).
Therefore, all $a_i$ are equal. Let $a_i = a$ for some $a \in k$.
This implies that $\Phi(\xi_i) = \xi_i + a \frac{\partial
F}{\partial \xi_j} = B_a(\xi_i)$ as desired. \hfill $\square$

\begin{proposition} \label{semi} Let $\gg$ be a finite dimensional simple Lie
superalgebra over $k$.

(i) Assume that $\gg$ is of type I and $\gg \neq \gp \gs \gl (2|2)$.
The map $\delta : \bG_m \to \bAut(\gg; \gg_0^{ss})
$ resulting from the standard $\Z$-grading of $\gg$ is an
isomorphism.

(ii) Assume that $\gg$ is of type II. The map  $\delta : \boldsymbol{\mu}_2 \to
\bAut(\gg; \gg_0^{ss})$ resulting from the standard
$\Z/2\Z$-grading of $\gg$ is an isomorphism.

(iii) If $\gg = \gp \gs \gl (2|2)$, then $\bAut(\gg;
\gg_0^{ss}) = \bSL_2^{\text{\rm out}}$.

(iv) If $\gg$ is of  Cartan type and $\gg \neq S'(2l)$, $\gg \neq
H(2l)$, the map $\delta : \bG_m \to \bAut(\gg;
\gg_0^{ss})$ resulting from the standard $\Z$-grading of
$\gg$ is an isomorphism.

(v) If $\gg = H(2l)$, then $\delta \times \beta: \bG_m \times \bG_a
\to \bAut(\gg; \gg_0^{ss})$ is an isomorphism.

(vi) Let $\gg = S'(2l)$. Then $\bAut(\gg; \gg_0^{ss})
\simeq \boldsymbol{\mu}_2$, where $\boldsymbol{\mu}_2$ corresponds
to the group $\Ad(\pm I) \subset \Aut_k (\gg)$.
\end{proposition}

\noindent {\bf Proof.}\footnote{The crucial ideas within this proof
are due to Serganova \cite{S}.}

(i) In this case $\gg_1$ and $\gg_{-1}$ are nonisomorphic
irreducible $\gg_0^{ss}$-modules.
 We must show that the closed embedding
$\delta : \bG_m \to \bAut(\gg; \gg_0^{ss})$ is surjective.
 Let $\varphi \in \bAut(\gg; \gg_0^{ss})$.
Then $\varphi_{| \gg_{\bar{1}}} : \gg_1 \oplus \gg_{-1} \to \gg_1
\oplus \gg_{-1}$ is a homomorphism of $\gg_0^{ss}$-modules
(because $\varphi$ fixes $\gg_0^{ss}$ pointwise) .

We now show that $\varphi(\gg_1) = \gg_1$. Consider the natural
projections $\pi_{\pm 1}: \gg_1 \oplus \gg_{-1} \to \gg_{\pm 1}$.
If $(\pi_{\pm 1} \circ \varphi) (\gg_1 \oplus \gg_{-1}) = 0$, then
$\varphi (\gg_{\pm 1}) \subseteq \gg_{\mp 1}$; which is impossible
given that $\varphi$ is injective and  $\gg_1$ and $\gg_{-1}$ are
nonisomorphic irreducible $\gs$-modules. We may therefore assume
that $(\pi_{ 1} \circ \varphi ) (\gg_1 \oplus \gg_{-1}) = \gg_{
1}$ and $(\pi_{- 1} \circ \varphi ) (\gg_1 \oplus \gg_{-1}) =
\gg_{- 1}$. If $\ker((\pi_{-1} \circ \varphi )_{|\gg_1}) = 0$,
then $\gg_1 \simeq \im((\pi_{-1} \circ \varphi )_{|\gg_1})
\subseteq \gg_{-1}$ As explained above, this forces $\gg_{-1}
\simeq \gg_1$ contrary to our assumption. Thus $\ker((\pi_{-1}
\circ \varphi )_{|\gg_1}) = \gg_1$ which implies $\varphi (\gg_1)
\subseteq \gg_1$. Similarly
 $\varphi (\gg_{-1}) = \gg_{-1}$.

By Schur's Lemma we conclude that ${\varphi}_{|\gg_{\pm 1}}=
\lambda_{\pm 1} \Id$ for some scalars $\lambda_{\pm 1} \in
k^{\times}$. But since $[\gg_1,\gg_{-1}] = \gg_0$, we obtain that
$\lambda_1 \lambda_{-1} = 1$ (because $\varphi$ fixes
$\gg_0^{ss}$), hence that $\varphi$ fixes $\gg_0$. This
completes the proof of (i).

(ii) The reasoning is similar to that of (i) above.

(iii) It is clear from the definition that $\bSL_2^{\text{\rm
out}} $ is a subgroup of $\bAut(\gg; \gg_0^{ss})$, so
(iii) comes down to showing that $\rho : \bSL_2\rightarrow
\bAut(\gg; \gg_0^{ss})$ is surjective. Let $\sigma \in
\bAut(\gg; \gg_0^{ss})$. Since we have a $\gg_0 $-module
isomorphism $\phi: \gg_1  \to \gg_{-1} $, we may apply
 Schur's Lemma to $\pi_{-1}\sigma i_{-1}$,
$\pi_{1}\sigma i_{1}$, $\phi^{-1} \pi_{-1}\sigma i_{1}$, and
$\phi^{-1} \pi_{1}\sigma i_{-1}$, where $i_{\pm 1}: \gg_{\pm 1}
\to \gg_{-1}  \oplus \gg_1$ and $\pi_{\pm 1}: \gg_{-1} \oplus
\gg_1 \to \gg_{\pm 1}$ are the natural inclusions and projections.
As a result we obtain that $\sigma \left(
\begin{array}{cc} A& B
\\ C & D  \end{array}
\right) = \left( \begin{array}{cc} A& \alpha B + \beta \psi (C) \\
\gamma \psi (B) + \delta C & D  \end{array} \right)$, for some
$\alpha, \beta, \gamma, \delta \in k $. Now using that
$\psi([B,C])= [\psi(B),\psi(C)]$ and
$$
\sigma \left( \left[ \left(
\begin{array}{cc} 0& B \\ 0 & 0  \end{array} \right), \left(
\begin{array}{cc} 0& 0 \\ C & 0  \end{array} \right) \right]
\right)= \left[ \sigma  \left( \begin{array}{cc} 0& B \\ 0 & 0
\end{array} \right) , \sigma \left( \begin{array}{cc} 0& 0 \\ C &
0  \end{array}  \right)  \right]
$$
we find $\alpha \delta - \beta
\gamma = 1$. Therefore $\sigma \in \bSL_2^{\text{\rm out}}$.

(iv) and (v) In these cases we use the fact that the standard
$\Z$-grading $\gg_{-1} \oplus \gg_0 \oplus ... \oplus \gg_s$ of
$\gg$ is such that $\gg_i$ is  an irreducible
$\gg_0^{ss}$-module  for $i=-1$ or $s$. To see  that the
morphisms are surjective we reason as follows. Let $\sigma \in
\bAut(\gg; \gg_0^{ss})$. Using Lemma \ref{nilp} we find
$\sigma_n \in \bN$ and $\sigma_0 \in \bAut(\gg)$, such that $\sigma_0(\gg_i)
= \gg_i$ and $\sigma = \sigma_0 \sigma_n$. Therefore
$\sigma_n (\gg_0^{ss}) = \sigma_0^{-1}\sigma
(\gg_0^{ss}) \subset \gg_0$ and by the definition of
$\bN$, $\sigma_{n|\gg_0^{ss}} = \Id$. By Lemma
\ref{cart} we see that $\sigma_n = \Id$ if $\gg \neq H(2l)$ and
$\sigma_n = \beta_a$ for some $a \in k$ if $\gg = H(2l)$.
Replacing $\sigma$ by $\sigma \sigma_n^{-1}$ if necessary we
conclude that $\sigma (\gg_i) = \gg_i$. Then by Schur's Lemma for
$j=-1, r$, $\sigma_{|\gg_{j}} = \lambda_j \Id$, for some
$\lambda_j \in k$. Now using $[\gg_{i+1}, \gg_{-1}] = \gg_i$ for
$i=r-1,r-2,...,-1$ we find $\sigma_{|\gg_{i}} = \lambda_i \Id$ and
$\lambda_i = \lambda_r (\lambda_{-1})^{r-i}$.

(vi) Let $\sigma \in \bAut(\gg; \gg_0^{ss})$. Again from
Lemma \ref{nilp} we have $\sigma = \sigma_0 \sigma_n$, where
$\sigma_n \in \bN$. By appealing to Lemma \ref{cart},
and reasoning as in the proof of (iv) above, we conclude that
$\sigma = \sigma_0$ . In the present case however $\gg$ does not
have a standard $\Z$-grading, but we get around this point by
making use of the explicit description of $\sigma_0$ given in
Table 1 of \cite{GP}. We have $\sigma_0 = \Ad (X)$ for some $X \in
\mbox{SL}_{2l}(k)$. Since $\sigma_{|\gg_0} = \Id$ then $X = \Id$
or $X = -\Id$. This completes the proof. \hfill $\square$

\subsection{The group $\bAut(\gg,\Pi_0$)}

This group measures the automorphisms of $\gg$ that induce
symmetries of the Dynkin diagram of $(\gg_0^{ss},
\gh_{\gg_0^{ss}})$ with respect to the chosen base $\Pi_0$
of $\Delta_{\gg_0^{ss}}$. By definition $\bAut(\gg ,\Pi_0):=
\{ \varphi \in \bAut(\gg)\; | \; \varphi(\gg_0^{ss}) =
\gg_0^{ss}, \varphi(\gh_{\gg_0^{ss}}) =
\gh_{\gg_0^{ss}}, \, \text{and} \, \varphi^* (\Pi_0)= \Pi_0
\}$.

\begin{proposition}\label{F}
Let $\bAut^0 (\gg)$ be the connected component of the identity of
$\bAut (\gg)$, and let $\F_k = \bAut (\gg)/\bAut^0 (\gg) $ be the
corresponding (finite constant)  group of connected components of
$\bAut (\gg)$ (see Theorem 4.1 in \cite{GP}).

(i) If $\sigma \in \bAut^0 (\gg)$ is  such that
$\sigma(\gg_0^{ss}) = \gg_0^{ss}$, then
$\sigma_{|\gg_0^{ss}}$ is in the connected component of the
identity  of $\bAut (\gg_0^{ss})$ (namely
 the Chevalley group of adjoint type of $\gg_0^{ss}$).

(ii) $\bAut(\gg, \Pi_0)  \cap \bAut^0(\gg) = \bH\big(\bAut (\gg ;
\gg_0^{ss}) \cap \bAut^0 (\gg)\big)$.

(iii) The restriction of the canonical map $j: \bAut(\gg) \to
\F_k$ to $\bAut(\gg, \Pi_0)$ is surjective.
\end{proposition}
\noindent {\bf Proof.} (i) We follow the descriptions of $\bAut^0 (\gg)$ and  $\G_{\bar{0}}^{sad}$
provided in Table 1 of \cite{GP} and the first table at the end of this paper, respectively. Let
 $\sigma \in \bAut^0(\gg)$ be such that $\sigma(\gg_0^{ss}) = \gg_0^{ss}$. We first consider the case when $\gg$ is of type I or II.
Then $\sigma$ is a product of three automorphisms: an element
$\sigma_0$ of $\G_{\bar{0}}^{sad}$; an element $j_k(\lambda)$ of
$\bG_m$, $\lambda \in k^{\times}$ (in fact, $j_k(\lambda) =
\delta_{\lambda}$);
 and (in the case $\gg = \gp \gs  \gl
(2|2)$) an element $\theta$ of $\bSL_2^{\text{\rm out}}$. But
$j_k(\lambda)_{|\gg_0^{ss}} = \theta_{|\gg_0^{ss}} = \Id$. On the
other hand, the restriction of $\G_{\bar{0}}^{sad}$ to
$\gg_0^{ss}$ is precisely the connected component of the identity
of $\bAut (\gg_0^{ss})$, from which the assertion follows. Let now
$\gg$ be of Cartan type. Then $\sigma = \sigma_{\lambda} \sigma_0
\sigma_n$, where (for $\gg = H(2l)$) $\sigma_{\lambda} \in \bG_m$
corresponds to the multiplication by $\lambda \in k^{\times}$ in
$\Lambda(2l)$ (in fact, $\sigma_{\lambda} = \delta_{\lambda}$),
$\sigma_0 \in \G_{\bar{0}}^{sad}$, and $\sigma_n \in \bN$. Then
since $\sigma$, $\sigma_{\lambda}$, and $\sigma_0$ leave invariant
$\gg_0^{ss}$, so it does $\sigma_n$. Thus $\sigma_{n|\gg_0^{ss}}=
\Id$, and, as before, we conclude that $\sigma_{|\gg_0^{ss}} =
\sigma_{0|\gg_0^{ss}}$ is in the connected component of the
identity of $\bAut (\gg_0^{ss})$.

(ii) We first check that $\bH$ is a subgroup of $\bAut (\gg,
\Pi_0)$. Let $\sigma \in \bH$. Since $\gh_{\gg_0^{ss}} \subset
\gg^0$ and $\sigma_{|\gg^0} = \Id$ we have that
$\sigma_{|\gh_{\gg_0^{ss}}}= \Id$. In particular, $\sigma^*(\Pi_0)
= \Pi_0$. Since $\gg_0^{ss} = \gh_{\gg_0^{ss}}\oplus
\left(\bigoplus_{\alpha_0 \in \Delta_{\gg_0^{ss}}}
(\gg_0^{ss})^{\alpha_0}\right)$ and $(\gg_0^{ss})^{\alpha_0} =
\gg^{\alpha} \cap \gg_0$ ( this last by Lemma \ref{compatible}) we
see that $\sigma$ acts as a multiplication by a constant
$\lambda(\alpha_0) \in k^{\times}$ on each
$(\gg_0^{ss})^{\alpha_0}$. In particular, $\sigma(\gg_0^{ss})
\subseteq \gg_0^{ss}$. This shows that $\bH \subset \bAut (\gg,
\Pi_0).$

Clearly $\bH$ normalizes $\bAut (\gg; \gg_0^{ss}) \cap \bAut^0
(\gg)$. Thus $\bH(\bAut (\gg; \gg_0^{ss}) \cap \bAut^0 (\gg))$
 is a closed subgroup of $\bAut (\gg)$. The inclusion
 $\bH(\bAut (\gg; \gg_0^{ss}) \cap \bAut^0 (\gg)) \subset \bAut(\gg , \Pi_0) \cap \bAut^0(\gg)$ is clear because $\bH$ is
 connected. Let $\sigma \in \bAut(\gg, \Pi_0)
\cap \bAut^0(\gg)$. By (i) we see that $\sigma_{|\gg_0^{ss}}$ is an
inner automorphisms of $ \gg_0^{ss}$ that stabilizes
$\gh_{\gg_0^{ss}}$ and $\Pi_0$ (this last via the $^*$ action). Thus
$\sigma_{|\gg_0^{ss}} = \Ad(x)_{|\gg_0^{ss}}$ for some $x \in \bT$,
where $\bT$ is the maximal torus of $\bG$ corresponding  to
$\gh_{\gg_0^{ss}}$. Note that $\Ad(\bT) \subset \bH$. Because $\bH$
is connected, we have $\bH \subset \bAut^0(\gg)$. Thus
$\Ad(x)^{-1}\sigma \in \bAut(\gg; \gg_0^{ss}) \cap \bAut^0(\gg)$ as
desired.

(iii) This statement follows by a case-by-case verification. For
example, for  $\gg = \gs \gl (m|n)$, we see that the
supertransposition $S : \left(
\begin{array}{cc} A& B \\ C & D  \end{array} \right)\mapsto \left(
\begin{array}{cc} -A^t& C^t \\ - B^t & -D^t
\end{array} \right)$ is in $\bAut(\gg , \Pi_0)$, and that
$j(S)$ generates $\F_k$. \hfill $\square$

\begin{corollary}
$\bAut(\gg , \Pi_0) / \bH (\bAut (\gg; \gg_0^{ss}) \cap \bAut^0
(\gg)) \simeq \F_k$.

\end{corollary}

\begin{remark}
{\rm We have $\bAut (\gg; \gg_0^{ss}) \subset  \bAut^0 (\gg)$ for
all $\gg \neq \gp \gs \gq (n)$. If $\gg = \gp \gs \gq (n)$, then
$\delta_{-1} \in \bAut (\gg; \gg_0^{ss})$ and $ \delta_{-1} \notin
 \bAut^0 (\gg)$.

This follows from a case-by-case verification using the first
table in the Appendix. If $\bAut (\gg; \gg_0^{ss})$ is connected
there is nothing to proof. Otherwise $\gg$ is of type II or $\gg =
S'(2l)$, and $\bAut (\gg; \gg_0^{ss}) \simeq \boldsymbol{\mu}_2$
is generated by the element $\delta_{-1}$ or $\Ad(-I)$. If $\gg
\neq \gp \gs \gq (n)$ then $\delta_{-1}$ is of the form $\Ad_{\gg}
(x)$ where $\Ad_{\gg}$ is the morphism of \cite{GP} used to
describe $\bAut^0 (\gg)$.\footnote{The morphism $\Ad$ defined in
Lemma \ref{ad}) can roughly be thought as the part of $\Ad{\gg}$
that arises from the semisimple part of $\gg$.} For example, if
$\gg = D(\alpha)$,  then $x = (I,I,-I)$. Furthermore, it is easy
to check that for $\gg = \gp \gs \gq (n)$, $\delta_{-1} =
\Ad_{\gg}(x)$ has no solutions for $x \in \bSL_n$.}
\end{remark}

\section{Abstract automorphisms of $\gg(R)$ and its universal
central extension} \setcounter{equation}{0}

\begin{lemma} \label{r-auto}
$\Aut_k(\gg (R))= \Aut_R(\gg (R)) \rtimes \Aut_k (R)$.

\end{lemma}
\noindent {\bf Proof.} Since $\gg$ is central and perfect (see
Proposition 7.1 in \cite{GP}), the Lemma follows from a
superversion of Lemma 4.4 in \cite{ABP} (see also Corollary 2.28
in \cite{BN}). More precisely, for a $k$-algebra automorphism
$\varphi : \gg (R) \to \gg (R) $ there exists a unique
$\widetilde{\varphi} \in \Aut_k (R)$ for which $\varphi(rx) =
\widetilde{\varphi}(r)\varphi(x)$ for all $x \in \gg $ and $r \in
R$. It is also evident that every element of $\Aut_k (R)$ lifts
naturally to and element of $\Aut_k (\gg(R))$. This leads to the
split exact sequence
$$
1 \to \Aut_R(\gg (R)) \to \Aut_k (\gg (R)) \to \Aut_k (R) \to 1.
$$ \hfill $\square$

\medskip

 Let  $\bG$ and
$\Ad$ be as in Lemma \ref{ad}. The (abstract) group $\Ad_R \,\bG(R)
\subset \bAut(\gg)(R) = \Aut_R (\gg(R))$ is in general much smaller
than the group of $R$-points of the the quotient group $\bG/{\rm
ker}(\Ad) \simeq \bG_{\bar{0}}^{\text{\rm sad}}$. While the group
$\Ad \,\bG(R)$ is quite explicit, $\bG_{\bar{0}}^{\text{\rm
sad}}(R)$ is not. The following Theorem shows that, for a large
class of objects in $\kalg$, an explicit and concrete description of
$\Aut_k (\gg(R))$ can still be achieved.

\begin{theorem}  \label{main} Assume  the object $R $ in $\kalg$ is a
Noetherian domain\footnote{Or more generally, that Spec$(R)$ is connected and
admits a rational point.} with trivial Picard group (for example
$R$ factorial).
 Then $\Aut_k (\gg(R))$ is generated by the subgroups $\Aut_k (R)$, $\Ad \bG(R)$,
$\bAut(\gg , \Pi_0)(R)$, $\bN(R)$, and $\bAut(\gg;
\gg_0^{ss})(R)$.
\end{theorem}

\noindent {\bf Proof.} Let $\sigma \in \Aut_k (\gg(R))$. By Lemma
\ref{r-auto}, we may assume that  $\sigma \in \Aut_R (\gg(R))$. Let
$\gg_{\bar{0}} = \gg_0 \oplus \gg_2 \oplus... \oplus \gg_{2r}$ be
the fixed $\Z$-grading of $\gg_{\bar{0}}$.  Using Lemma \ref{nilp},
we can write $\sigma$ as a product $\sigma_0 \sigma_n$, where
$\sigma_0(\gg_0(R)) = \gg_0(R)$ and $\sigma_n \in \bN(R)$. Replacing
$\sigma$  by $ \sigma \sigma_n^{-1}$, we may then assume that
$\sigma(\gg_{0}(R))\subseteq \gg_{0}(R)$. Since $\gg_0^{ss} =
[\gg_0, \gg_0]$, we conclude that $\sigma(\gg_0^{ss}(R))\subseteq
\gg_0^{ss}(R)$. Given our assumptions on $R$, and by taking Lemma
\ref{ad} into consideration, we can appeal to the conjugacy theorem
of regular maximal abelian $k$-diagonalizable subalgebras of
$\gg_0^{ss}(R)$ (\cite{P2} Theorem 1(ii)(a)) for the existence of an
element of $\Ad \bG(R)$ taking $\sigma(\gh_{\gg_0^{ss}})$ to
$\gh_{\gg_0^{ss}}$. We may thus assume that $\sigma$ stabilizes
$\gh_{\gg_0^{ss}}$. Lemma \ref{star} implies that the contragradient
automorphism $\sigma^*$ of $\gh_{\gg_0^{ss}}^*$ stabilizes
$\Delta_{\gg_0^{ss}}$.  By means of the Weyl group of $(\gg_0^{ss},
\gh_{\gg_0^{ss}})$, whose elements we can recreate as restrictions
to $\gh_{\gg_0^{ss}}$ of elements of $\Ad \bG(k)$, we may further
assume that $\sigma^* (\Pi_0) = \Pi_0$. Let ${\frak m}$ be a maximal
ideal of $R$ for which $R/{\frak m} \simeq k$. Then $\sigma \otimes
1 \in \Aut_R ((\gg \otimes R) \otimes_R R/{\frak m}) \simeq \Aut_k
(\gg).$ Clearly $\sigma \otimes 1$, when viewed
 as an element of  $\bAut(\gg)$, is in
fact an element of $\bAut(\gg , \Pi_0)$. Let $\widetilde{\sigma}$
denote the $R$-linear extension of this element to $\Aut_R (\gg
(R))$. Then, after replacing $\sigma$ by $\widetilde{\sigma}^{-1}
\sigma$, we may assume that $\sigma$ fixes
 $\gh_{\gg_0^{ss}}$ pointwise, hence that
 $\sigma$
stabilizes $\gg^{\alpha} \otimes R$ for every $\alpha$ in
$\Delta_{\gg_0^{ss}}$ (see Lemma \ref{star}).

We now proceed by a case-by-case reasoning using Lemmas \ref{t2}
and \ref{t3}.

{\it Case 1:  $\gg = \gp \gs \gl (2|2)$; or $\gg_0^{ss}$ is simple, $\gg \neq H(2k)$.} In
this case we use Lemma \ref{t2} (i) and fix a basis $B$ of $Q_{\gg}$
such that $p(B) \cap \Delta_{\gg_0^{ss}}$ is a base of
$\Delta_{\gg_0^{ss}}$. By multiplying $\sigma$ with an
element of $\bH(R)\subset \bAut(\gg , \Pi_0)(R)$ (see Lemma \ref{F} (ii)) we may assume that $\sigma$ fixes a  set of
generators $e_{\alpha}$ of  $(\gg_0^{ss})^{p(\alpha)} = \gg_0 \cap \gg^{\alpha}$ (see Lemma \ref{t2} (iii)),
 for any $\alpha \in B$. Since $\sigma$ is $R$-linear,
it fixes $\gg_0^{ss}(R)$, and thus is in $\bAut(\gg;
\gg_0^{ss})(R)$.

{\it Case 2: $\gg_0^{ss}$ is not simple and $\gg \neq \gp
\gs \gl (2|2)$.} Now we use Lemma \ref{t3}. For our chosen base
$\Pi = \{ \alpha_1, \alpha_2,...,\alpha_l\}$ of $\Delta$, we fix
$e_i \in \gg^{\alpha_i}$, $f_i \in \gg^{-\alpha_i}$, and
$\alpha_i^{\vee} = [e_i, f_i].$ Since the spaces $\gg^{\alpha_i}$ are
$1$-dimensional, after
   multiplying by an element of $\bH(R)$
 we may assume that $\sigma$ fixes $e_i$, $f_i$, and $\alpha_i^{\vee}$.
Since these generate $\gg$ and $\sigma$ is $R$ linear, we have
$\sigma = \Id$.

{\it Case 3: $\gg = H(2l)$, $l \geq 3$}. In this case using the
description of $\Delta$ provided in the Appendix we see that
$\gg_{-1} \cap \gg^{\alpha} \neq 0$ iff $\alpha  = \pm
\varepsilon_i$, $i = 1,...,l$. Moreover, the spaces $\gg_{-1} \cap
\gg^{\pm \varepsilon_i}  = \gg_{-1}^{p(\pm \varepsilon_i)}$ are
one-dimensional. Multiplying $\sigma$ by an element of $\bH(R)$ we
may assume that
 $\sigma_{|\gg_{-1} \cap \gg^{\varepsilon_i}}=\Id$  for every
$i=1,...,l$. Let $r \in R^{\times}$ be such that
$\sigma_{|\gg_{-1} \cap \gg^{-\varepsilon_1}}=r\Id$ (strictly
speaking $\sigma_{|\gg_{-1}(R) \cap \gg^{-
\varepsilon_1}(R)}=r\Id$). From $[\gg_0,\gg_{-1}]=\gg_{-1}$ and
the fact that the spaces $\gg_{i} \cap \gg^{\alpha}$, $i = -1,0$,
are at most one-dimensional, we easily conclude that
$\sigma_{|\gg_{-1} \cap \gg^{-\varepsilon_i}}=\sigma_{|\gg_0 \cap
\gg^{-\varepsilon_i-\varepsilon_j}}=r\Id$, $\sigma_{|\gg_0 \cap
\gg^{\varepsilon_i-\varepsilon_j}}= \Id$, and $\sigma_{|\gg_0 \cap
\gg^{\varepsilon_i+\varepsilon_j}}=r^{-1}\Id$ for every $1 \leq i
\neq j \leq l$. This completely determines $\sigma \in \Aut_R
(\gg(R))$ since every automorphism of $\gg$ is uniquely determined
by its restriction on $\gg_{-1}$ (see \cite{S}). In order to
explicitly express $\sigma$ we apply the change of coordinates $
\eta_i:=\frac{1}{\sqrt{2}}(\xi_i + \sqrt{-1}\xi_{i+l});\;
\eta_{i+l}:=\frac{1}{\sqrt{2}}(\xi_i - \sqrt{-1}\xi_{i+l})$. Then
$D_{\eta_{i+l}}  \in \gg_{-1} \cap \gg^{\varepsilon_i}$ and
$D_{\eta_{i}}  \in \gg_{-1} \cap \gg^{\varepsilon_{i+l}}$. In
terms of the new coordinates we have that $\sigma(D) =
\bar{\sigma}D\bar{\sigma}^{-1}$, where $\bar{\sigma}$ is the
linear automorphism of $\Lambda(2l)(R)$ given by the matrix
$A_{\sigma} = \left(\begin{matrix}  r^{-1}I & 0 \\ 0 &
I\end{matrix} \right)$, i.e. $\bar{\sigma}(\eta) = A_{\sigma}
\eta$. It then follows that $\sigma \in \bAut(\gg , \Pi_0)(R)$.
\hfill $\square$

\medskip
\begin{remark} {\rm Let $r \in R^\times,$ and consider the quadratic
extension $\widetilde{R}=R[r^{\frac{1}{2}}]$ of $R$. Then
$\widetilde{R}/R$ is finite \'etale (in fact Galois). Let
$\widetilde{\sigma} = \Ad \left(
\begin{matrix}  r^{\frac{1}{2}} I & 0 \\ 0 & r^{-\frac{1}{2}} I\end{matrix} \right)
\delta_{r^{\frac{1}{2}}}$. Then $\widetilde{\sigma} \in \bAut(\gg)
(\widetilde{R})$ is such that $\widetilde{\sigma}$ stabilizes
$\gg(R)$ and $\widetilde{\sigma}_{|\gg(R)} = \sigma$, where $\sigma$
is the automorphism of $\gg(R)$ determined by the matrix
$A_{\sigma}$ in Case 4 above. According to Proposition \ref{F} we
should be able to write $\sigma$ as an $R$-point of the product of
the groups $\bH$ and $\bAut (\gg; \gg_0^{ss})$. The $R$-points of
the product group are in general a larger group than the naive
product of the $R$-points of the respective groups. In fact, $\Ad
\left( \begin{matrix}  r^{\frac{1}{2}} I & 0 \\ 0 & r^{-
\frac{1}{2}} I\end{matrix} \right) \in \bH (\widetilde{R})$ and
$\delta_{r^{\frac{1}{2}}} \in \bAut (\gg,
\gg_0^{ss})(\widetilde{R})$, which shows that $\sigma$ is an
$R$-point of the product group.}

\end{remark}
\medskip

\begin{remark} {\rm Let $\widetilde{\gg(R)}$ be the universal central
extension of the $k$-Lie superalgebra $\gg (R)$. A result of Neher
(Corollary 2.8 in \cite{Ne}) implies that $\Aut_k (\gg(R)) =
\Aut_k (\widetilde{\gg(R)})$. Theorem \ref{main} can also be applied
to this last group. }

\end{remark}

The following corollary (of the proof) of Theorem \ref{main}, is
useful for the representation theory of $\gg(R)$ (notably in the
case of $R=k[t, t^{-1}]$, which corresponds to the untwisted
affine Kac-Moody superalgebras). It provides a description of the
subgroup of the $k$-automorphisms of $\gg(R)$ that leave invariant
the category of weight $(\gg(R),
\gh_{\gg_0^{ss}})$-modules, i.e., all $\gg(R)$-modules $M$
for which which $M = \bigoplus_{\lambda \in
\gh^*_{\gg_0^{ss}}}M^{\lambda}$.

\begin{corollary}Let $R$ be as in Theorem \ref{main}. Let $\mbox{\bf T}$ be the maximal
torus of $\bG$ corresponding to $\gh_{\gg_{0}^{ss}}.$ The subgroup
of $\Aut_k (\gg(R))$ consisting of all automorphisms $\sigma$ for
which $\sigma(\gh_{\gg_0^{ss}}) = \gh_{\gg_0^{ss}}$, is generated by
$\Ad (N_{\bG}(\mbox{\bf T}))(k)$ together with $\Aut_k (R)$,
$\bAut(\gg , \Pi_0)(R)$, and $\bAut(\gg; \gg_0^{ss})(R)$.
\end{corollary}
\noindent {\bf Proof.}
 We may assume that $\sigma \in \Aut_R
(\gg(R))$. The group $\Ad(N_{\bG}(\mbox{\bf T}))$ accounts for the
action of the Weyl group used in the proof of Theorem \ref{main}.
Now if $\sigma$ stabilizes $\gh_{\gg_0^{ss}}$, we reason as in the
proof of Theorem \ref{main} to conclude that $\sigma$ belongs to
the subgroup of $\Aut_R\gg(R)$ which is generated by the subgroups
prescribed by the Corollary. \hfill $\square$

\newpage

\centerline{ {\bf Tables}}
\bigskip

\noindent {\bf 1}. The groups $\bG_{\bar{0}}^{sad}$ and $\bAut
(\gg; \gg_0^{ss})$.

$$
\begin{tabular}{|c|c|c|c|}
\hline

${\gg}$ & $\bG_{\bar{0}}^{sad} $ &  $\bAut (\gg;
\gg_0^{ss})$  \\ \hline $\gs \gl (m|n)$ & $(\mbox{\bf SL}_m
\times \mbox{\bf SL}_n
)/ (\boldsymbol{\mu}_m \times \boldsymbol{\mu}_n)$ &  $\mbox{\bf G}_m $  \\
\hline

$\gp \gs \gl (n|n), n>2$ & $(\mbox{\bf SL}_n \times \mbox{\bf SL}_n )/
(\boldsymbol{\mu}_n \times \boldsymbol{\mu}_n)$ & $\mbox{\bf G}_m $\\
\hline

$\gp \gs \gl (2|2)$ & $(\mbox{\bf SL}_2  \times \mbox{\bf SL}_2)/\boldsymbol{\mu}_2$ &
$\mbox{\bf SL}_2$
\\ \hline

$\gs \pp (n)$ & $\mbox{\bf SL}_n / \boldsymbol{\mu}_n $ & $\mbox{\bf G}_m $  \\
\hline

$\gp\gs \gq (n)$ & $\mbox{\bf SL}_n  / \boldsymbol{\mu}_n $ & $\boldsymbol{\mu}_2 $  \\
\hline

$\go \gs \gp (2l|2n), l\neq 1$ & $(\mbox{\bf SO}_{2l} \times
\mbox{\bf Sp}_{2n}) / \boldsymbol{\mu}_2 $ &
$ \boldsymbol{\mu}_2 $\\
\hline

$\go \gs \gp (2|2n)$ & $(\mbox{\bf SO}_{2} \times \mbox{\bf
Sp}_{2n}) / \boldsymbol{\mu}_2 $ &
$\mbox{\bf G}_m $\\
\hline

 $\go \gs \gp (2l+1|2n)$ & $\mbox{\bf SO}_{2l+1} \times
\mbox{\bf Sp}_{2n}$ & $\boldsymbol{\mu}_2 $\\ \hline

$F(4)$ & $(\mbox{\bf Spin}_{7} \times \mbox{\bf SL}_{2})/\boldsymbol{\mu}_2$ &$\boldsymbol{\mu}_2 $\\
\hline

$G(3)$ & $\mbox{\bf G}_{2} \times \mbox{\bf SL}_{2}$ & $\boldsymbol{\mu}_2 $ \\
\hline

$D(\alpha)$  & $ (\mbox{\bf SL}_{2} \times \mbox{\bf SL}_{2}
\times\mbox{\bf SL}_{2})/ (\boldsymbol{\mu}_2 \times
\boldsymbol{\mu}_2)$ & $\boldsymbol{\mu}_2 $ \\ \hline

$W(n)$ & $\mbox{\bf GL}_n$& $\mbox{\bf G}_m  $\\
\hline

$S(n)$ & $\mbox{\bf SL}_n$& $\mbox{\bf G}_m  $ \\
\hline

$S'(2l)$ & $\mbox{\bf SL}_{2l}$& $ \boldsymbol{\mu}_2  $ \\
\hline

$H(2l)$ & $\mbox{\bf SO}_{2l}/\boldsymbol{\mu}_2$& $\mbox{\bf G}_a \times \mbox{\bf G}_m  $ \\
\hline

$H(2l+1)$ &$\mbox{\bf SO}_{2l+1}$&  $\mbox{\bf G}_m$ \\
\hline
\end{tabular}
$$

\bigskip

\bigskip

\noindent {\bf 2}. The groups $\bN$.
\bigskip

Two alternative definitions of the unipotent groups $\bN$
 have been provided in \S 3.  The table below gives
an explicit description of $\bN$ in terms of the $n$-dimensional
vector space $V:=k\xi_1 \oplus... \oplus k \xi_n$, where
$\xi_1,...,\xi_n$ are odd variables, i.e. $\xi_i^2 = 0, \xi_i
\xi_j = - \xi_j \xi_i$ if $i \neq j$. Recall that the  additive
affine group of a finite dimensional $k$-space $U$ is denoted by
$U_a$, that is $U_a = \Hom_k(S(U^*), - )$. For the proofs we refer
the reader to Lemmas 6.1 and 6.2 in \cite{GP}.

$$
\begin{tabular}{|c|c|}

\hline

$\gg$ & $\bN$ \\

\hline
$W(n)$ & $\mbox{\bf Hom}(V, \oplus_{i\geq 1}\Lambda^{2i}V)$ \\
\hline

$S(n)$ & $(\oplus_{i\geq 1} (V^* \otimes \Lambda^{2i+1}V)/\Lambda^{2i}V)_a$ \\
\hline

$S'(n), n=2l$ & $(\oplus_{i\geq 1} (V^* \otimes \Lambda^{2i+1}V)/\Lambda^{2i}V)_a$ \\
\hline

$H(n)$ & $(\oplus_{i\geq 2} \Lambda^{2i}V)_a$ \\
 \hline

\end{tabular}
$$

\newpage
\appendix{\noindent {\bf Appendix: Cartan subsuperalgebras and root
systems of the Cartan type Lie superalgebras}}
\bigskip

We first recall some generalities about the Cartan type Lie
superalgebras. By $W(n)$ we denote the (super)derivations of the
Grassmann algebra $\Lambda(n) :=\Lambda(\xi_1,...,\xi_n)$ over
$k$. Every element $D$ of $W(n)$ is of the form $D = \sum_{i=1}^n
P_i(\xi_1,...,\xi_n) \frac{\partial}{\partial \xi_i}$ where by
definition $\frac{\partial}{\partial \xi_i} (\xi_j) =
\delta_{ij}$. Both $\Lambda(n)$ and $W(n)$ have natural gradings
$\Lambda(n) = \oplus_{i=0}^n \Lambda(n)_i$ and  $W(n) =
\oplus_{j=-1}^{n-1} W(n)_j$, where $\Lambda(n)_i:= \{
P(\xi_1,...,\xi_n) \; | \; \deg P = i\}$ and
$W(n)_{j}:=\{\sum_{i=1}^n P_i \frac{\partial}{\partial \xi_i} \; |
\; \deg P_i = j+1\}$. For any Lie subsuperalgebra $\gs$ of $W(n)$
we set $\gs_j:=\gs \cap W(n)_j$.

In this appendix we will use the following explicit description of
the Cartan type Lie subsuperalgebras  $S(n)$, $S'(n)$, $H(n)$, and
$\widetilde{H}(n)$, of $W(n)$:
\begin{eqnarray*}
S(n) & = & \Span_k \left\{ \frac{\partial f}{\partial \xi_i}
\frac{\partial}{\partial \xi_j} +  \frac{\partial f}{\partial
\xi_j} \frac{\partial}{\partial \xi_j}\; | \; f \in \Lambda(n),
i,j = 1,...,n\right\},\\
 S'(n)& = & \Span_k \left\{ (1 - \xi_1 ...
\xi_n)\left(\frac{\partial f}{\partial \xi_i}
\frac{\partial}{\partial \xi_j} + \frac{\partial f}{\partial
\xi_j} \frac{\partial}{\partial
\xi_j}\right)\; | \; f \in \Lambda(n), i,j = 1,...,n\right\},\\
\widetilde{H}(n)& = & \Span_k \left\{D_f:=\sum_i \frac{\partial
f}{\partial \xi_i} \frac{\partial}{\partial \xi_i} \; | \; f \in
\Lambda(n), f(0)=0, i,j = 1,...,n\right\},\\
\widetilde{H}(n) & = & H(n) \oplus kD_{\xi_1...\xi_n}.
\end{eqnarray*}
We have also that $[D_f,D_g] = D_{\{f,g\}}$ where
$\{f,g\}:=(-1)^{\deg f}\sum_{i=1}^n \frac{\partial f}{\partial
\xi_i} \frac{\partial g}{\partial \xi_j}$.

In what follows we give an explicit description of specific
Cartan subsuperalgebras of $W(n), S(n), S'(n)$, and $H(n)$. This
description can serve as a complement to the root structures
provided in Appendix A in \cite{Pen}.

\medskip \noindent
{\it Case 1: $\gg = W(n)$, $n \geq 2$}. In this case the elements
$h_i = \xi_i \frac{\partial}{\partial \xi_i}$ form a basis for a
Cartan subsuperalgebra $\gh$ of $\gg$. In particular, $\gh =
\gh_{\bar{0}}$ is a subalgebra of $\gg_0 \simeq \gg \gl (n)$. The
elements in $\gh_{\bar{0}}^*$ that form the dual basis to $h_i$
will be denoted by $\varepsilon_i$. The root system of $\gg$ is
$$
\Delta = \{ \varepsilon_{i_1} +...+ \varepsilon_{i_k},
\varepsilon_{i_1} +...+  \varepsilon_{i_k}- \varepsilon_{j} \; |
\; i_r \neq i_s, j \neq i_r, 0 \leq k \leq n - 1, 1 \leq j \leq
n\}
$$

\medskip \noindent
{\it Case 2: $\gg = S(n), S'(n)$ ($n = 2l$ in the second case), $n
\geq 3$}. Now we consider the Cartan subsuperalgebra $\gh$ spanned
by $h_i - h_j = \xi_i \frac{\partial}{\partial \xi_i} - \xi_j
\frac{\partial}{\partial \xi_j}$. We have again that $\gh =
\gh_{\bar{0}}$ is a subalgebra of $\gg_0 \simeq \gs \gl (n)$.
Denote by $\varepsilon_i$ the images in $\gh_{\bar{0}}^*$ of the
basis dual to $h_i$, $\varepsilon_1 +...+\varepsilon_n = 0$.

The root system of $\gg$ is
$$
\Delta = \{ \varepsilon_{i_1} +...+ \varepsilon_{i_k},
\varepsilon_{i_1} +...+  \varepsilon_{i_l}- \varepsilon_{j} \; |
\; i_r \neq i_s, j \neq i_r, 1 \leq k \leq n - 2, 0 \leq l \leq
n-1, 1 \leq j \leq n\}
$$

\medskip \noindent
{\it Case 3: $\gg = H(2l), l \geq 3$}. We fix $\gh$ to be the
Cartan subsuperalgebra of $\gg$ spanned by the elements
$$
\{ D_{\xi_{i_1}...\xi_{i_r} \xi_{i_1+l}...\xi_{i_r+l}}\; | \; 1
\leq r \leq l-1, 1\leq i_1< ...< i_r \leq l\}.
$$

In this case we have $\gh = \gh_{\bar{0}}$ and $\gh_0^{ss}:= \gh \cap
\gg_0 = \Span_k \{ D_{\xi_{i} \xi_{i+l}}\; | \; 1\leq i \leq l \}$
is a Cartan subalgebra of $\gg_0 \simeq \gs \go (2l)$. Moreover, $\gh_{\bar{0}} = \gh_{\gg_{0}^{ss}} \oplus \gh^2$,
where $\gh^2 \subset \gg^2 = \oplus_{i\geq 1}\gg_{2i}$. Let
$h_i:=\sqrt{-1}D_{\xi_i \xi_{i+l}}$ and set $\varepsilon_i$ to be
the basis of $(\gh_0^{ss})^*$ dual to $h_i$, i.e. $\varepsilon_i (h_j) =
\delta_{ij}$ (note that the choice of $\varepsilon_i$
is meaningful because $\alpha_{|\gh^2} = 0$ for every $\alpha \in \Delta$). The root system $\Delta = \Delta_{\bar{0}}\cup
\Delta_{\bar{1}}$ of $\gg$ is described by
$$
\Delta = \{\varepsilon_{i_1} +...+  \varepsilon_{i_t}-
\varepsilon_{j_1} -...-  \varepsilon_{j_s}\; | \; i_r \neq i_p,
j_r \neq j_p, i_r \neq j_p, 0 \leq t,s \leq l\}.
$$

\medskip \noindent
{\it Case 4: $\gg = H(2l+1), l \geq 2$}. We fix $\gh$ to be the
Cartan subsuperalgebra of $\gg$ spanned by the elements
$$
\{ D_{\xi_{i_1}...\xi_{i_r} \xi_{i_1+l}...\xi_{i_r+l}\xi_{2l+1}},
D_{\xi_{i_1}...\xi_{i_s} \xi_{i_1+l}...\xi_{i_s+l}}\; | \; 0 \leq
r \leq l -1, 1 \leq s \leq l, 1\leq i_1< ...< i_r \leq l\}.
$$

Now we have $\gh = \gh_{\bar{0}} \oplus \gh_{\bar{1}}$ where
$\gh_{\bar{0}}$ is spanned by the elements
$D_{\xi_{i_1}...\xi_{i_s} \xi_{i_1+l}...\xi_{i_s+l}}$, while
$\gh_{\bar{1}}$ is spanned by the elements
$D_{\xi_{i_1}...\xi_{i_s} \xi_{i_1+l}...\xi_{i_s+l}\xi_{2l+1}}$. In particular,
$\gh_0^{ss}:= \gh \cap \gg_0 = \Span_k \{ D_{\xi_{i} \xi_{i+l}}\; | \;
1\leq i \leq l \}$ is a Cartan subalgebra of $\gg_0 \simeq \gs \go
(2l+1)$. We have again that $\gh_{\bar{0}} = \gh_{\gg_{0}^{ss}} \oplus \gh^2$,
where $\gh^2 \subset \gg^2 = \oplus_{i\geq 1}\gg_{2i}$. As in Case 3 we set $h_i:=\sqrt{-1}D_{\xi_i \xi_{i+l}}$ and
$\varepsilon_i \in (\gh_0^{ss})^*$ with $\varepsilon_i (h_j) =
\delta_{ij}$. Then the root system of $\gg$ is
$$
\Delta = \{\varepsilon_{i_1} +...+  \varepsilon_{i_t}-
\varepsilon_{j_1} -...-  \varepsilon_{j_s}\; | \; i_r \neq i_p,
j_r \neq j_p, i_r \neq j_p, 0 \leq t,s \leq l\}.
$$

In what follows we describe the graded root spaces $\gg^{\alpha}
\cap \gg_i$ of $\gg = H(2l)$ and $\gg = H(2l+1)$. For this
description it is convenient to use the following coordinate
change: $ \eta_i:=\frac{1}{\sqrt{2}}(\xi_i +
\sqrt{-1}\xi_{i+l}),\; \eta_{i+l}:=\frac{1}{\sqrt{2}}(\xi_i -
\sqrt{-1}\xi_{i+l})$, $\eta_{2l+1}:=\xi_{2l+1}$ (the last in the
case $\gg = H(2l+1)$). Using the new coordinates we have that if
$\gg = H(2l+\epsilon)$, $\epsilon = 0$ or $1$, then $D_f=
\sum_{i=1}^l \frac{\partial f}{\partial \eta_{i+l}}
\frac{\partial}{\partial \eta_i} + \sum_{i=1}^l \frac{\partial
f}{\partial \eta_{i}}\frac{\partial}{\partial \eta_{i+l}} +
\epsilon \frac{\partial f}{\partial
\eta_{2l+1}}\frac{\partial}{\partial \eta_{2l+1}}$, and $\{f,g\} =
(-1)^{\deg f}\left(\sum_{i=1}^l \frac{\partial f}{\partial
\eta_{i+l}} \frac{\partial g}{\partial \eta_i} + \sum_{i=1}^l
\frac{\partial f}{\partial \eta_{i}}\frac{\partial g}{\partial
\eta_{i+l}} + \epsilon\frac{\partial f}{\partial
\eta_{2l+1}}\frac{\partial g}{\partial \eta_{2l+1}}\right)$.

For $I=(i_1,...,i_t)$, $1 \leq i_1 <...<i_t \leq l$, we set
$\widehat{I}:=(i_1+ l,...,i_t+l)$, $|I|:= i_1+...+i_t$,
$\varepsilon_I:= \varepsilon_{i_1}+ ... +\varepsilon_{i_t}$, and
$\eta_I:=\eta_{i_1}...\eta_{i_t}$. Fix now $I$ and $J$ such that
$I \neq J$. For $\gg = H(2l)$ we have that $D_{\eta_I
\eta_{\widehat{J}}} \in \gg^{\varepsilon_I - \varepsilon_J} \cap
\gg_{|I| +|J| - |I \cap J| - 1}$. If $I \cap J = \emptyset$, then
the set
$${\mathcal
B}_{I,J}:=\left\{ D_{\eta_I \eta_K \eta_{\widehat{J}}
\eta_{\widehat{K}} } \; | \; K\cap I = \emptyset, K \cap J =
\emptyset \right\}
$$
forms a basis of $\gg^{\varepsilon_I -
\varepsilon_J}$ and the set
$${\mathcal B}_{I,J, i}:={\mathcal
B}_{I,J} \cap \gg^i = \left\{ D_{\eta_I \eta_K \eta_{\widehat{J}}
\eta_{\widehat{K}} } \; | \; K\cap I = \emptyset, K \cap J =
\emptyset, |K|= \frac{i-1-|I|-|J|}{2}\right\} $$ forms a basis of
$\gg^{\varepsilon_I - \varepsilon_J}\cap \gg_i$ (we have
$\gg^{\varepsilon_I - \varepsilon_J}\cap \gg_i \neq 0$ iff
$i-1-|I|-|J|$ is even).

In the case of $\gg = H(2l+1)$ and $I \cap J = \emptyset$, the set
${\mathcal B}_{I,J} \cup{\mathcal B}'_{I,J}$ forms a basis of
$\gg^{\varepsilon_I - \varepsilon_J}$, where
$$
{\mathcal B}'_{I,J}:= \left\{ D_{\eta_I \eta_K \eta_{\widehat{J}}
\eta_{\widehat{K}}\eta_{2l+1}} \; | \; K\cap I = \emptyset, K \cap
J = \emptyset \right\}.
$$
Furthermore, $\gg^{\varepsilon_I - \varepsilon_J}\cap \gg_i$ has
${\mathcal B}_{I,J, i}$ as a basis if $i-1-|I|-|J|$ is even and
$${\mathcal B}'_{I,J, i}:={\mathcal B}'_{I,J} \cap \gg_i = \left\{ D_{\eta_I \eta_K
\eta_{\widehat{J}} \eta_{\widehat{K}}\eta_{2l+1} } \; | \; K\cap I
= \emptyset, K \cap J = \emptyset, |K|=
\frac{i-|I|-|J|}{2}\right\}$$ as a basis if $i-1-|I|-|J|$ is odd.

Dimitar Grantcharov

Department of Mathematics

San Jose State University

San Jose, CA 95192-0103

USA

E-mail: grantcharov@math.sjsu.edu

\bigskip

Arturo Pianzola

Department of Mathematical and Statistical Sciences

University of Alberta

Edmonton, Alberta T6G 2G1

CANADA

E-mail: a.pianzola@math.ualberta.ca


\begin{thebibliography}{99}



\bibitem[ABP]{ABP} B. Allison, S. Berman, A. Pianzola, Covering algebras II:
Isomorphisms of Loop algebras, {\it J. reine angew. Math.} {\bf
571} (2004), 39--71.

\bibitem[BN]{BN} G. Benkart and E. Neher, The centroid of Extended Affine and Root
Graded Lie algebras, {\it J. Pure and Appl. Algebra} {\bf 205} (2006), 117--145.

\bibitem[B]{Borel} A. Borel, Linear Algebraic Groups, Benjamin, New
York (1969).

\bibitem[DG]{DG} M. Demazure, P. Gabriel, Groupes alg\'ebriques.
Tome I: G\'eom\'etrie alg\'ebrique, g\'en\'eralit\'es, groupes
commutatifs.  North-Holland Publishing Co., Amsterdam, 1970.


\bibitem[GP]{GP} D. Grantcharov, A. Pianzola, Automorphisms and twisted loop
algebras of finite dimensional simple Lie superalgebras, {\it Int.
Math. Res. Not.} {\bf 73} (2004), 3937--3962.

\bibitem[GiPi1]{GiPi1} P. Gille and A. Pianzola, Isotriviality of
torsors over Laurent polynomials rings,  C. R. Acad. Sci. Paris,
Ser. I {\bf 340} (2005) 725--729.

\bibitem[GiPi2]{GiPi2} P. Gille and A. Pianzola,  Galois cohomology and forms of algebras over Laurent polynomial rings,
  Math. Annalen (in press).

\bibitem[K]{K} V. Kac, Lie superalgebras, {\it Adv. Math.} {\bf 26} (1977), 8--96.

\bibitem[Ne]{Ne} E. Neher, An introduction to universal central
extensions of Lie superalgberas, {\it Math. Appl.} {\bf 555}
(2003) 141--166.

\bibitem[Pen]{Pen} I. Penkov, Characters of strongly generic irreducible Lie superalgebra representations, {\it Internat. J. Math.} {\bf 9} (1998), 331--366.

\bibitem[PS]{PS} I. Penkov, V. Serganova, Generic irreducible representations of finite-dimensional
Lie superalgebras, {\it Internat. J. Math.} {\bf 5} (1994), 389--419.


\bibitem[P1]{P1} A. Pianzola, Automorphisms of toroidal Lie algebras and their central quotients, {\it J. Algebra Appl.} {\bf 1} (2002), 113--121.

\bibitem[P2]{P2} A. Pianzola, Locally trivial principal homogeneous
spaces and conjugacy theorems for Lie algebras, {\it Journal of
Algebra} {\bf 275} (2004) 600--614.

\bibitem[P3]{P3} A. Pianzola, Line bundles and conjugacy theorems
for toroidal Lie algebras, {\it C. R. Acad. Sci. Canada} {\bf 22}
(2000) 125-128.

\bibitem[P4]{P4} A. Pianzola, Affine Kac-Moody Lie algebras as
torsors over the punctured line, Indagationes Mathematicae N.S.
{\bf 13}(2) (2002) 249--257.

\bibitem[P5]{P5} A. Pianzola, Vanishing of $H^1$ for Dedekind rings and applications to loop
algebras, C. R. Acad. Sci. Paris, Ser. I {\bf 340} (2005),
633--638.


\bibitem[Sch]{Sch} M. Scheunert, Invariant supersymmetric multilinear forms and the
Casimir elements of $P$-type Lie superalgebras, {\it J. Math.
Phys.} {\bf 28} (1987), 1180-1191.


\bibitem[S]{S}  V. Serganova, Automorphisms of simple Lie superalgebras, {\it Math. USSR Izvestiya} {\bf 24} (1985), 539-551.



\end{thebibliography}
\end{document}